\iftrue %

\documentclass[12pt,a4paper]{article}
\usepackage{time,psfig,amsmath,amssymb} %
\usepackage[T1]{fontenc}
\usepackage{mathptmx}
\def\today{\number\year\space
 \ifcase\month\or January\or February\or March\or April\or May\or June\or
   July\or August\or September\or October\or November\or December\fi
     \space\number\day
}
\title{\vspace{-15mm}\bf Some explicit\\ badly approximable pairs}
\author{Keith Briggs,
BTexact Technologies\\
Adastral Park, Antares 2 pp5\\
Suffolk IP5 3RE, UK\\
{\tt Keith.Briggs@bt.com}
}
\date{\today} %

\begin{document}
\maketitle
\begin{abstract}
\noindent
I consider the Diophantine approximation problem of sup-norm
simultaneous rational approximation with common denominator of a pair
of irrational numbers, and
compute explicitly some pairs with large 
approximation constant.    One of these pairs is the most badly
approximable pair yet computed.
\end{abstract}

\else %

\documentclass[12pt,doublespacing]{elsart}
\usepackage{time,psfig,amsmath,amssymb}
\usepackage[T1]{fontenc}
\usepackage{mathptmx}
\journal{Journal of Number Theory}
\begin{document}
\begin{frontmatter}
\title{Some explicit badly approximable pairs}
\author{Keith Briggs}
\address{BTexact Technologies,
Adastral Park, Antares 2 pp5,
Suffolk IP5 3RE, UK\\
{\tt Keith.Briggs@bt.com}
}

\begin{abstract}
I consider the Diophantine approximation problem of sup-norm
simultaneous rational approximation with common denominator of a pair
of irrational numbers, and
compute explicitly some pairs with large 
approximation constant.    One of these pairs is the most badly
approximable pair yet computed.
\end{abstract}
\begin{keyword}
Simultaneous Diophantine approximation \sep badly approximable
\end{keyword}
\end{frontmatter}

\fi

\def\Z{\mathbb Z}
\def\R{\mathbb R}
\def\Q{\mathbb Q}
\medmuskip=1mu plus 0.5mu minus 1.5mu  %

\noindent
The theory of approximation of a single irrational number by rationals is well
known, and for our purposes the relevant facts may be summarized as follows.
We measure the goodness of approximation of the rational number $p/q$ to
$\alpha$ by $c(\alpha,p,q)\equiv q|q\alpha-p|$. For each irrational $\alpha$
(without loss of generality, we may assume $0\!<\!\alpha\!<\!1$) we know by
Dirichlet's thereom that there are infinitely many rationals $p/q$ such that
$|\alpha-p/q|\!<\!1/q^2$, or $c(\alpha,p,q)\!<\!1$.  It is therefore of
interest to ask how small one may make $\gamma$ in $c(\alpha,p,q)\!<\!\gamma$
before this property fails to hold.  The {\sl approximation constant\/} of
$\alpha$ is thus defined as $c(\alpha)\equiv\liminf_{q\rightarrow\infty}
c(\alpha,p,q)$.  Here, of course, for each $q$ we choose the $p$ whch minimizes
$c(\alpha,p,q)$.  Numbers $\alpha$ with a large $c(\alpha)$ are hard to
approximate by rationals.  The {\sl one-dimensional Diophantine approximation
constant\/}, defined as $c_1=\sup_{\alpha\in\R}c(\alpha)$, has the value
$1/\sqrt{5}$, attained at $\alpha=(\sqrt{5}-1)/2$.

Otherwise expressed, this means that $c_1$ is the unique number such that for
each $\epsilon\!>\!0$, the inequality $c(\alpha,p,q)<c_1+\epsilon$ has
infinitely many rational solutions $p/q$ for all $\alpha$, whereas there is at
least one $\alpha$ such that $c(\alpha,p,q)<c_1-\epsilon$ has only finitely
many rational solutions.

These results completely solve the problem of rational approximation in one
dimension, but by contrast the situation in two or more dimensions is much more
complex and in fact the value of the analogous constant $c_n$ for $n\!\geqslant
\!2$ is unknown \cite{Szekeres84a,Szekeres85a}.

We wish to simultaneously approximate a pair of irrationals by a pair of
rationals with common denominator and to measure the closeness of approximation
by the maximum error in the two components, so we make the definitions:
for ${p}=(p_1,p_2)\in \Z^2, q\in\Z, \alpha=(\alpha_1,\alpha_2)\in\R^2$, let
$$
c(\alpha,p,q)=q\,\max\left(\,|q\alpha_1-p_1|^2,|q\alpha_2-p_2|^2\,\right)
$$
and
$$
c({\alpha})=\liminf_{q\rightarrow\infty}\;\{\;c({\alpha},p,q),\; {p}\in \Z^2,\; q\in \Z\;\}.
$$
The two-dimensional (sup-norm) simultaneous Diophantine approximation 
constant is then
$$
c_2=\sup_{\alpha\in \R^2}\;c(\alpha).
$$
Despite much work over the last few decades
\cite{Cassels55a,Adams69a,Adams80a,Cusick74b,Cusick83a,Szekeres84a,Szekeres85a},
the value of $c_2$ is unknown, though folk-lore suggests that its value is
$2/7$.    Adams \cite{Adams69a} has shown that this is the correct value if we
restrict the pair $(\alpha_1,\alpha_2)$ to cubic number fields, but his result
does not give us a constructive procedure to identify pairs with large
$c(\alpha)$.

Here, however, I use a theorem of Cusick together with high-precision numerical
computation to explicitly compute examples of such pairs.   These have
potential applications to numerical simulation studies of dynamical systems on
the 2-torus, where $(\alpha_1,\alpha_2)$ represent the winding number of
periodic orbits.

Cusick's construction makes use of the cubic number field $\Q(\theta)$, where
$\theta=2\cos(2\pi/7)$, of smallest positive discriminant, namely 49.  For
details on cubic fields and their integral bases, I refer to
\cite{Cohen93}.

The theorem of Cusick \cite{Cusick74b} states that for any integral basis
$\{1,\alpha,\beta\}$ of $\Q(\theta)$, we have $c^{*}\!\!<\!2/7$, where
$c^{*}$ is the infimum of those $c$ such that 
$$\left|\,x+\alpha y+\beta z\,\right|\;\max\,(y^2,z^2) \,<\, c$$ 
(with $y$ and $z$ not both zero) has infinitely many solutions in integers
$x,y,z$.  Additionally, for any $\epsilon\!>\!0$ there is an integral basis
$\{1,\alpha,\beta\}$ such that $$2/7-c^{*}(\alpha,\beta)\,<\,\epsilon$$ iff
\begin{itemize}
\item [1:] The continued fraction of $\theta$ has patterns 
$[\ldots,n_1,1,1,n_2,\ldots]$ with $n_1,n_2$ arbitrarily large; or,
\item [2:] The continued fraction of $\theta$ has patterns 
$[\ldots,n_1,2,n_2,\ldots]$ with $n_1,n_2$ arbitrarily large.
\end{itemize}
It is not known whether either of the last two conditions are satisfied.  Note
that this theorem relates to the dual problem to simultaneous Diophantine
approximation, namely approximation to zero by linear forms. Hence, it is not
immediately apparent that the upper bound of $2/7$ that it gives for $c^{*}$ it
defines is relevant to the problem of determining $c_2$.  However, from another
paper by Cusick (\cite{Cusick72a}, Corollary 1 on page 187), we have that for
the particular field $\Q(\theta)$, $c^{*}(\alpha,\beta)=c(\alpha,\beta)$ for all
integral bases.  Also, by a theorem of Davenport \cite{Davenport52a}, we have
$\sup c^{*}(\alpha,\beta)=\sup c(\alpha,\beta)$, where the $\sup$s are over all
irrational pairs, not necessarily in a cubic field.

Thus, if the above patterns in the continued fraction of $\theta$ do in fact
exist, Cusick's theorem gives us a way of finding explicit pairs (which together
with 1 form an integral basis of $\Q(\theta)$) with a value of $c$ close to
$2/7$.  Even if $n_1, n_2$ do not become arbitrarily large, just the presence
of some large values gives us potential candidates for very badly approximable
pairs.

From results in \cite{Cusick74b}, it follows that for an integral basis of
the form $\{1,p\theta+q\theta^2,r\theta+s\theta^2\},
\left(\begin{smallmatrix}p q\\ r s\end{smallmatrix}\right)\in PSL(2,\Z)$,
where $-q/p$ and $-s/r$ are rational approximants to $\theta$ obtained
by truncating the continued fraction at the points where condition 1 or 2
is satisfied, $c^{*}$ is explicitly given by
\begin{eqnarray}
c^{*}&=&1/\max\,\{\,|A+B+C|,|A-B+C|,|C-B^2/(4A)|,|A-B^2/(4C)|\,\}\nonumber\\
   &=&1/\max\,\{\,|A+B+C|,|A-B+C|,49/|4A|,49/|4C|\,\},\nonumber
\end{eqnarray}
where
\begin{eqnarray}
\left[\begin{array}{c}A\\B\\C\end{array}\right]
&=&
\left[\begin{array}{lll}
s^2&-rs&r^2\\
-2qs&ps+qr&-2pr\\
q^2&-pq&p^2\\
\end{array}\right]
\left[\begin{array}{c}a\\b\\c\end{array}\right]
\nonumber
\end{eqnarray}
with
\begin{eqnarray}
a_{\hphantom{1}}&=&(\theta_2^2-\theta^2)(\theta^2-\theta_1^2)\nonumber\\
b_{\hphantom{1}}&=&(\theta_2^2-\theta^2)(\theta_1-\theta)+(\theta_2-\theta)(\theta_1^2-\theta^2)\nonumber\\
c_{\hphantom{1}}&=&(\theta-\theta_2)(\theta_1-\theta)\nonumber\\
\theta_{\hphantom{1}}&=&2\cos(2\pi/7)\nonumber\\
\theta_1&=&2\cos(4\pi/7)\nonumber\\
\theta_2&=&2\cos(6\pi/7)\nonumber.
\end{eqnarray}

With this background, I can now state the main result of this paper: 
I have exactly computed over 2 million %
partial quotients of the continued fraction of $\theta$ (directly from
the defining cubic $x^3+x^2-2x-1$), and the required 
patterns do indeed occur, though very infrequently.   
The largest values of $c^{*}$, with the corresponding fractional parts of
$\alpha=p\theta+q\theta^2$ and $\beta=r\theta+s\theta^2$ occur at:

\begin{itemize}
\item [(A)] positions 57-60: $[\ldots,60,1,1,50,\ldots],\qquad c^{*}\approx 0.2851877$ %

\hspace{-5mm}\begin{minipage}{\linewidth}
{\small
$\alpha\approx$ 0.4563286858107963651609830446124431560745665647128596153008802
\\
$\beta \approx$ 0.4781573193903170892895817415258772866671562381178937772663665
}
\end{minipage}

\item [(B)] positions 2927-2930: $[\ldots,22,1,1,22,\ldots],\qquad c^{*}\approx 0.2853154$ %

\hspace{-5mm}\begin{minipage}{\linewidth}
{\small
$\alpha\approx$ 0.1554011929520066325796747316744656830061413509133865038820677 
\\
$\beta \approx$ 0.6003679362632065361061389158735863615694126556922931077332356
}
\end{minipage}
\item [(C)] positions 3626-3629: $[\ldots,272,1,1,215,\ldots],\qquad c^{*}\approx 0.2855726$ %

\hspace{-5mm}\begin{minipage}{\linewidth}
{\small
$\alpha\approx$ 0.6530646111210617321254054547968773238346090082060701183776580
\\
$\beta \approx$ 0.9410463762107594592302548739412493098027738320829952592216557
}
\end{minipage}
\item [(D)] positions 33877-33880: $[\ldots,81,1,1,78,\ldots],\qquad c^{*}\approx 0.2856261$ %

\hspace{-5mm}\begin{minipage}{\linewidth}
{\small
$\alpha\approx$ 0.9319638477108390366188499907354642637920661848031694636081724
\\
$\beta \approx$ 0.7032571495109702868148790086182835032528572663181375225766851
}
\end{minipage}
\item [(E)] positions 215987-215990: $[\ldots,124,1,1,129,\ldots],\qquad c^{*}\approx 0.2856678$ %

\hspace{-5mm}\begin{minipage}{\linewidth}
{\small
$\alpha\approx$ 0.4375520476578757564544576313180510209212270982522655674846137
\\
$\beta \approx$ 0.5646614639128419094417646922292433724548272488193131214134926
}
\end{minipage}
\item [(F)] positions 957740-957743 $[\ldots,460,1,1,415,\ldots],\qquad c^{*}\approx 0.28568046$ %

\hspace{-5mm}\begin{minipage}{\linewidth}
{\small
$\alpha\approx$ 0.6134980317071692745070006892224661159462079954445253478668675
\\
$\beta \approx$ 0.9411544329571988683307282702558980820407618535628393885417987
}
\end{minipage}
\item [(G)] positions 1650050-1650053: $[\ldots,648,1,1,666,\ldots],\qquad c^{*}\approx 0.2857082$ %

\hspace{-5mm}\begin{minipage}{\linewidth}
{\small
$\alpha\approx$ 0.4848739572889332951989678247806190621159456336657613155291560
\\
$\beta \approx$ 0.5404925035004667478257428539575752367424111926723566428410541
}
\end{minipage}
\end{itemize}

These calculations involve extremely large integer and floating-point numbers;
in case (G) the absolute values of the integers $p,q,r,s$ are of the order
$2^{3\times 10^6}$, and the calculation of $c^{*}$ requires floating-point
operations of about twice this precision.  In fact, these examples all come
from cases of Cusick's first condition, and $c^*$ is given by $49/|4A|$ or
$49/|4C|$.  Of course, the approximate decimal values for $\alpha,\beta$ given
above are insufficient to represent the true values, but these may be
reconstructed if required from the continued fraction of $\theta$.

\begin{figure}[tph]
  \begin{center}
    \leavevmode 
    \hbox{%
      \psfig{file=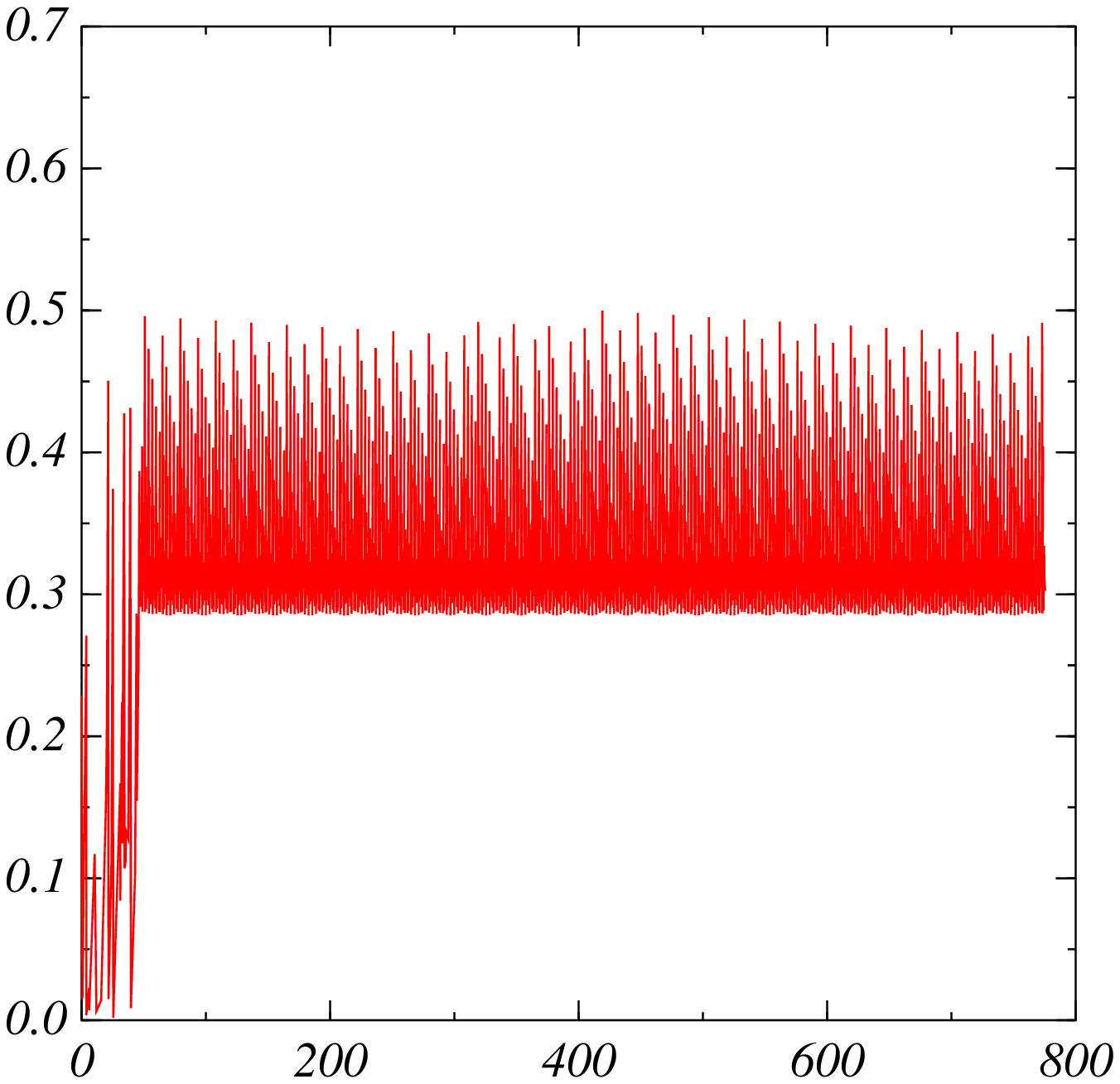,angle=0,width=9.5cm}
    }
    \hbox{%
      \psfig{file=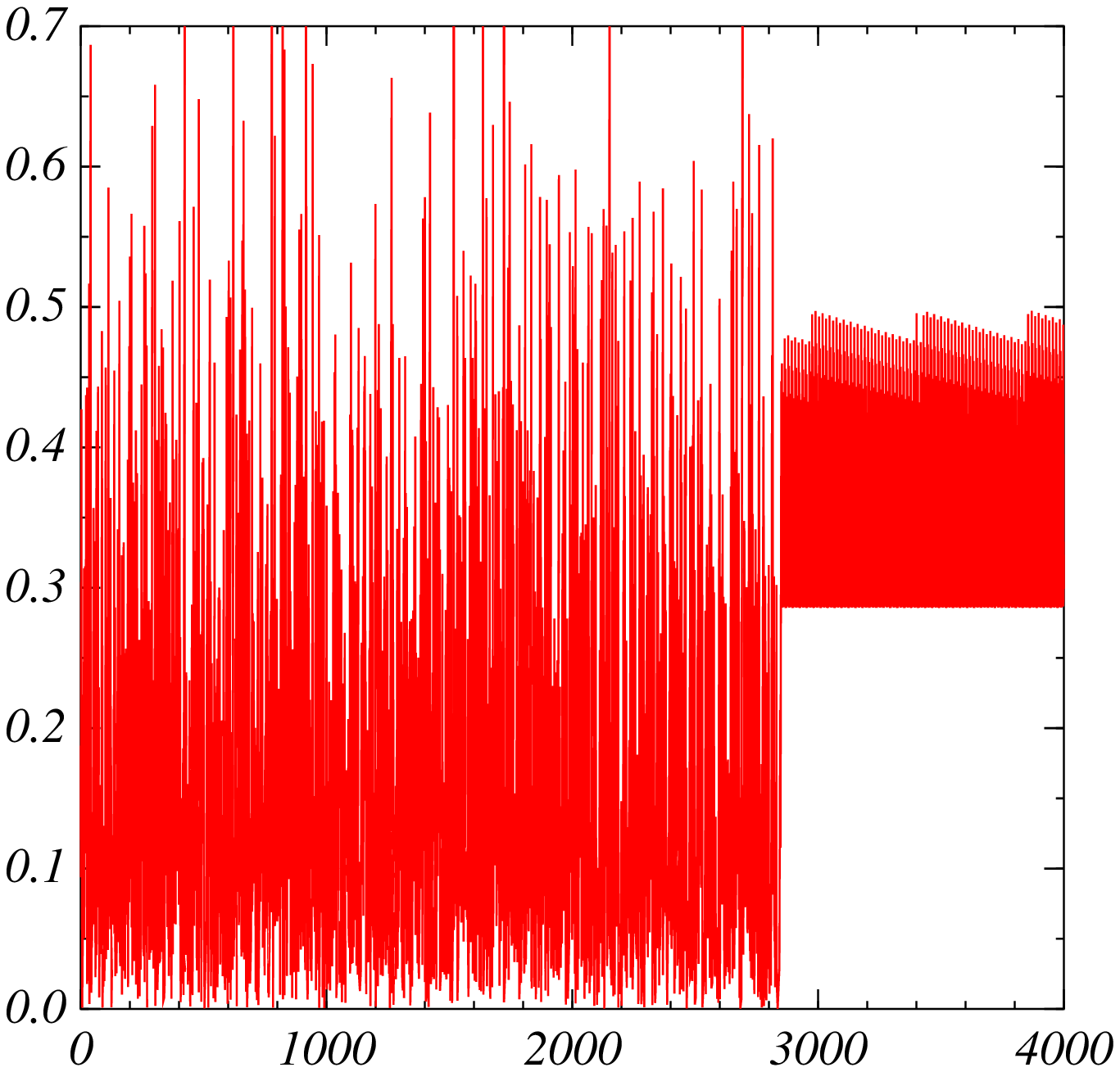,angle=0,width=9.5cm}
    }
    \caption{
      $c(\alpha,\beta,q)$ vs.\ $\log_{10}(q)$ at best approximants
      for two integral bases $(1,\alpha,\beta)$ of the field $\Q(\theta)$.
      Above: case (A), below: case (D). 
    }
    \label{fig4}
  \end{center}
\end{figure}

An independent verification of these results may be obtained by giving the
values $\alpha,\beta$ as input to a simultaneous Diophantine approximation
algorithm.    Such an algorithm finds all best simultaneous approximants up to
a given denominator.   For the computation of sup-norm best approximants,
an algorithm has been given by Furtw\"angler \cite{furtw28,briggs01a}.
Figure~\ref{fig4} shows the behaviour of the Furtw\"angler algorithm applied
the pairs (A) and (C) above. The approximation constant estimated from the
minimum $c$ after ignoring the initial transient is about $0.2856$, verifying
the more precise value of $c^{*}$ above.  But the chief point to be noted is
the extremely long initial transient.  Until a sufficient large denominator $q$
is reached, these pairs would in fact appear to be {\sl not\/} badly
approximable.

I have thus exhibited some explicit pairs which are very badly approximable by
rationals.  I believe that the value $0.2857082$ above is the largest
explicitly computed lower bound for the two-dimensional simultaneous
Diophantine approximation constant $c_2$.

The question remains open as to whether there are pairs 
(necessarily unrelated to the field $\Q(\theta)$)
with approximation constant larger than $2/7$.

\clearpage

\bibliography{furtw}
\vfill
\noindent{\tiny tantalum:home/kbriggs/Cusick/badly\_approx\_pairs.tex}
\hfill
\noindent{\tiny Typeset in \LaTeXe}
\end{document}